\documentclass{article}
\usepackage{amsmath}
\usepackage{tikz}
\usepackage{graphicx}
\title{Determination of the Representative Sample Size in Linear Regression}
\author{Anatoly Rayev}
\date{}
\begin{document}
\maketitle
\begin{center}
\textbf{ABSTRACT}
\end{center}

Very often, accuracy of analysis and forecasting (multiple coefficient of regression and residual means) obtained for a sample used to formulate a regression model is not equal to the accuracy achieved for another homogeneous sample. Indeed, accuracy of analysis and forecasting based on another sample is much worse. This is explained by the discrepancy between a postulated and a real model. This discrepancy is caused by the redundancy of the number of terms of an approximating series. It describes a phenomenon under study with the sample noise. To filter this noise it is necessary to correctly choose the terms of an approximating series and to determine their number, which in turn depends on sample size. It is possible to include in this approximating series not only independent variables but also various functions of the value of the independent variables (e.g., squares, cubes, algorithms, etc.). This paper gives the procedure for selection of an approximating series.

\vspace{1em}

\section*{1. PROBLEM}

The least squares estimate in regression analysis is:
\begin{equation}
b = (\bar{X} * \bar{X}')^{-1} * \bar{X}' * \bar{Y} \label{eq:1} \tag{1}
\end{equation}

of $\beta$ for the model $E(Y) = \bar{X} * \beta$ is an unbiased estimator only if the postulated model is the true model to consider (Draper and Smith 1981, p. 117), where

\begin{equation}
\bar{X} = \begin{bmatrix} 
1 & X_{11} & X_{12} & \dots & X_{1k} \\ 
\vdots & \vdots & \vdots & \vdots & \vdots \\ 
1 & X_{N1} & X_{N2} & \dots & X_{Nk} 
\end{bmatrix} 
\qquad 
\bar{Y} = \begin{bmatrix} 
Y_1 \\ 
Y_2 \\ 
\vdots \\ 
Y_N 
\end{bmatrix} \label{eq:2} \tag{2}
\end{equation}

N - is the number of observations \\
K - is the number of independent variables \\
$Y_i$ - is the value of a dependent variable in observation i \\
$X_{ij}$ - is the value of an independent variable j in observation i \\
$\bar{X}'$ - is the transposed matrix of $\bar{X}$ \\
$(\bar{X} * \bar{X}')^{-1}$ - is the inverse matrix of $(\bar{X} * \bar{X}')$ \\

If in the postulated model the number of independent variables exceeds the number of those in a true model, the coefficients of redundant independent variables (bj) are zero and the postulated model is consistent with a true model. In this case b) are also unbiased estimates. However, this is true only if the sample size is sufficiently large; in real life the sample size is always small.

To select the terms (independent variables) of an approximating series, we often use the "forward selection procedure." The largest possible number of terms is determined by the F-test, which gives the significance of residuals after the next term has been included in the approximating series.

If the sample size is small, an approximating series does not include all independent variables. In other words, we truncate a series. In this case a postulated model does not correspond to a real one, since we cannot include all relevant independent variables, and all estimates, including the F-test, are biased. In turn, an approximating series now includes redundant terms. In other words, if the value of F were not biased, we would select a smaller number of relevant terms. But since this value is biased, an approximating series includes more terms than are necessary. This redundancy means that the analysis of the interrelationship we are trying to describe includes noise which is imposed on the dependent variable. Very often this occurs when the sign of the coefficient of correlation between some independent and dependent variable is different from the sign of a regression coefficient (bj) in the regression equation.

When estimates are biased, the terms of an approximating series selected on the basis of a forward selection procedure can include unimportant independent variables, whereas terms of great significance can be excluded. Thus, the ranking of an approximating series of terms based on their influence on the dependent variable is disturbed. The selection of "garbage" terms begins when the F-test becomes biased. However, it is impossible to ascertain the exact moment when it occurs. That is why we can be confident only about the first term, which has the greatest correlation coefficient with a dependent variable Y.

Besides, an estimate of the coefficient of multiple determination and values of residuals does not correspond to the determination coefficient and values of residuals which were obtained when the regression equation was used for forecasting or analysis based on new (other) homogeneous samples. In this case the determination coefficient is lower and the residuals values are higher. In other words, errors increase. Thus, if a sample is small, it is impossible to use the parametric methods of estimation of regression coefficients and the value of the F-test. It is important to remember that the size of necessary (sufficient) sample depends on the number of independent variables and on the level of noise which is imposed on both the dependent and the independent variables.

Unfortunately, a priori noise values cannot generally be ascertained. This is why it is impossible to determine a priori the necessary sample size or the necessary number of terms of an approximating series. However, it is not possible to determine the number of terms to which the approximate series was cut either because all other dependent variables are not influencing or because the sample size is small. Hence, it is necessary to use non-parametric estimation methods.

The following simple example can illustrate the point. Assume we have one dependent variable, one independent variable and only three observations. We included only one independent variable because in this case it would be possible to use two-dimensional graphs. Thus, $Y_i$ and $X_i$ , $i = 1,3$ (see Figure 1).

Now we will formulate regression equations for two models:

\begin{equation}
YL = a_o + a_1 * X \label{eq:3} \tag{3}
\end{equation}

and

\begin{equation}
YQ = a_o + a_1 * X + a_2 * X^2 \label{eq:4} \tag{4}
\end{equation}

The first equation is a linear model; the second one is a quadratic function.

\begin{center}
\begin{tikzpicture}
\draw[->] (0,0) -- (4,0) node[right] {Y};
\draw[->] (0,0) -- (0,3) node[above] {X};
\node at (1, 1) {*};
\node at (1.5, 2.5) {*};
\node at (2.5, 1.5) {*};
\end{tikzpicture} \\
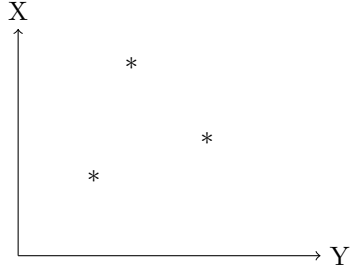
Figure 1. Values for Three Observations (Points).
\end{center}

which we do not have but which should be analyzed and predicted, will be approximately plotted in Figure 4.

\begin{center}
\begin{tikzpicture}
\draw[->] (0,0) -- (4,0) node[right] {Y};
\draw[->] (0,0) -- (0,3) node[above] {X};
\draw (0.2, 0.5) -- (3.8, 2.8);
\draw[thick] (0.5,0.5) parabola bend (2,2.8) (3,0.5);
\node at (1, 1) {\tiny *};
\node at (1.2, 1.6) {\tiny *};
\node at (1.5, 2.2) {\tiny *};
\node at (1.8, 1.2) {\tiny *};
\node at (2.1, 1.9) {\tiny *};
\node at (2.5, 2.6) {\tiny *};
\node at (2.8, 1.8) {\tiny *};
\node at (2.9, 2.1) {\tiny *};
\end{tikzpicture} \\
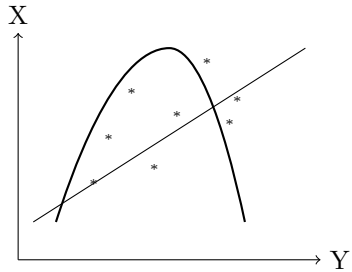
Figure 4. Observations of Linear and Quadratic Regression Models.
\end{center}

Hence, it is easy to see that the error for quadratic regression would be significantly higher in comparison to the linear model. The same situation will exist when we will try to formulate linear regression equations with more than one independent variable. To avoid this error, it is possible to use the F criterion, which estimates the significance of residual differences and truncates the series, adding the new members of the approximating series when this difference becomes insignificant.

However, as we have already emphasized, all estimates of the regression equation are unbiased only when the postulated model coincides with the true one. But this is impossible to know. That is why some variables are not included in the model or may have different relationships. But even if we know the true relationship, because of the limited sample size we cannot include all variables in the approximating series. Hence, the value of the F criterion will be biased and, as a result, the approximating series will include redundant members. That is why we get incorrect regression models with all inherent errors.

\section*{2. THE SOLUTION}

\subsection*{2.1 Review of Subsistence Methods}
Below we will formulate the two non-parametric algorithms that can be used to determine the sufficiency of sample size. Where the sample size is deemed to be sufficient, we can use the standard programs of regression analysis (for example, SPSS or SAS). Where the sample size is determined to be insufficient, these algorithms can be used to establish the permissible number of terms of an approximating series, which ensures that the regression estimates are unbiased.

\subsection*{Algorithm 1: "Data Splitting"}
If the sample size is sufficiently large, it can be divided into two parts. The first part is used to estimate the regression equation. The second part is used to test the quality of the estimation and is usually smaller in size. It is used to calculate residual mean squares D(k) or the coefficient of multiple determination, R(k), for the different number k of terms of the approximating series. As a rule, the selection of these terms is based on forward selection procedure:

\begin{equation}
D(k) = \frac{\sum_{m=1}^{M} (Y_m - Y_m(k))^2}{M} \label{eq:5} \tag{5}
\end{equation}

where:\\
$Y_m$ - the value of the dependent variable in observation m; \\
$Y_m(k)$ - the calculated value of the dependent variable of estimated regression function for observation m when the number of terms of the approximating series equals k; \\
$M$ - the number of observations in the second part of the sample.

Function D(k) can have one of three forms:

\begin{center}
\begin{tikzpicture}[scale=0.8]
\begin{scope}[xshift=0cm]
\draw[->] (0,0) -- (2.5,0) node[right] {K};
\draw[->] (0,0) -- (0,2.5) node[above] {D(k)};
\draw (0, 1.5) node[left] {S1 +};
\draw (1, 0) node[below] {K1} node {$+$};
\draw[thick] (0.2, 2.3) -- (1, 0.5) -- (1.8, 2.3);
\node at (1, -1) {(a)};
\end{scope}

\begin{scope}[xshift=4cm]
\draw[->] (0,0) -- (2.5,0) node[right] {K};
\draw[->] (0,0) -- (0,2.5) node[above] {D(k)};
\draw (0, 1.5) node[left] {S1 +};
\draw (1, 0) node[below] {K1} node {$+$};
\draw[thick] (0.2, 2.3) -- (1, 0.5) -- (1.5, 0.5) -- (2, 2.3);
\node at (1, -1) {(b)};
\end{scope}

\begin{scope}[xshift=8cm]
\draw[->] (0,0) -- (2.5,0) node[right] {K};
\draw[->] (0,0) -- (0,2.5) node[above] {D(k)};
\draw (0, 1.5) node[left] {S1 +};
\draw (1, 0) node[below] {K1} node {$+$};
\draw[thick] (0.2, 2.3) -- (1, 0.5) -- (2.2, 0.5);
\node at (1, -1) {(c)};
\end{scope}
\end{tikzpicture} \\
\vspace{1em}
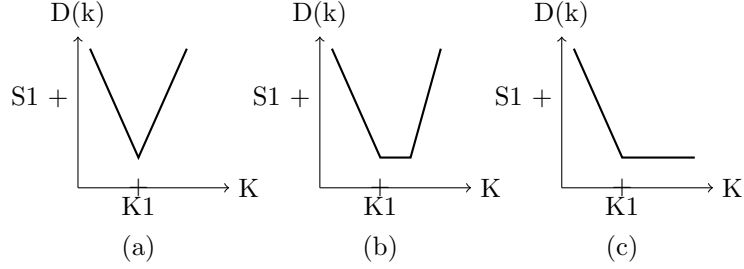
Figure 5. Possibility of Form Function D(k). Using resting data not used to calculate coefficients of regression, Figure 5(a) shows that the sample size is not enough to include all significant terms in the approximating series. Figures 5(b) and 5(c) show that the sample size is sufficient. In all three cases the number of terms of the approximate series equals (K1) and the unbiased estimate of residual is (S1).
\end{center}

Figure 5a. shows that the sample size is insufficient to include all significant terms (factors) in the approximating series. Figures 5b. and 5c. show that the sample size is sufficient. In all three cases the number of terms of the approximate series is K1 and the unbiased estimate of residual variance is S1.

If the calculation of D(k) is based on the sample which is also used to estimate regression function, then D(k) will have the following form:

\begin{center}
\begin{tikzpicture}
\draw[->] (0,0) -- (3,0) node[right] {K};
\draw[->] (0,0) -- (0,2.5) node[above] {D(k)};
\draw[thick] (0.2,2.3) to[out=280,in=170] (2.8,0.2);
\end{tikzpicture} \\
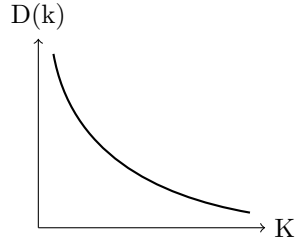
Figure 6. Dependence Residuals from Number of Terms of Approximating Series, Calculated on This Same Sample. This is the same data that was used to calculate coefficient of regression. The curve does not have a minimum. Consequently, it is impossible to determine the number of the approximate series (K1) and estimate the residual (S1).
\end{center}

It is clear that in this case it is impossible to determine K1 and to estimate the residual variance S1. If we construct the curve describing the relationship between the coefficient of multiple correlation R(k) and the number of terms of the approximating series using the second part of the sample, then the graph of this function will have the following form:

\begin{center}
\begin{tikzpicture}
\draw[->] (0,0) -- (3,0) node[right] {K};
\draw[->] (0,0) -- (0,2.5) node[above] {R(k)};
\draw (0,1.5) node[left] {1 +} -- (2.5,1.5);
\draw (1.5,0) node[below] {K1} node {$+$};
\draw[thick] (0.2,0.2) to[out=70,in=180] (1.5,1.3) to[out=0,in=140] (2.8,0.3);
\end{tikzpicture} \\
\end{center}

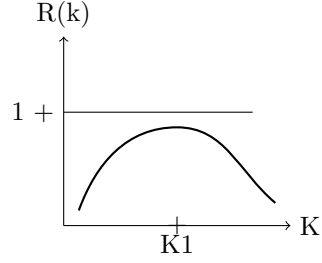
Figure 7. Dependence Multiple Coefficients Correlation from Terms of Approximating Series, Calculated on Different Sample. Using testing data not used to calculate coefficients of regression, the curve has a maximum and the number of terms of the approximate series is (K1).

If this function is calculated on the basis of the sample used to determine the regression coefficients, then the graph has the following form:

\begin{center}
\begin{tikzpicture}
\draw[->] (0,0) -- (3,0) node[right] {K};
\draw[->] (0,0) -- (0,2.5) node[above] {R(k)};
\draw (0,2) node[left] {1 +} -- (3,2);
\draw (0,0) node[left] {0};
\draw[thick] (0.2,0.2) to[out=60,in=185] (2.8,1.9);
\end{tikzpicture} \\
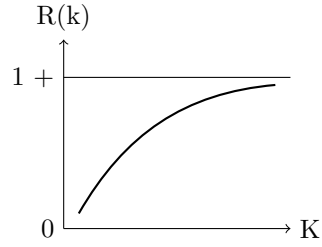
Figure 8. Dependence Multiple Coefficients Correlation from Terms of Approximating Series, Calculated on This Same Sample. Using the same data that was used to calculate coefficients of regression, the curve does not have a maximum. Consequently it is impossible to determine the number of the approximate series (K1).
\end{center}

It is clear that it is impossible to determine the value of K1. In practice, it is not always feasible to have the additional sample used for testing. Sometimes to do so would involve additional costs. That is why it is important to obtain unbiased estimates of the regression function on the basis of the same sample.

\subsection*{Algorithm 2 "Technique of Cross-validation" (Snee 1977, pp 415-418)}
Suppose we have a sample of size N (i=1, N). Estimate the function Y(k) using a sample of size N-1. We leave one observation to test and calculate the value of Di(k) (see formula 2) for this sample point (i=1). Of course, the significance of D(k) for one observation is zero. Assume, then, that we return this point (observation) back to the sample and choose the second point for testing. That is, we estimate the function Y(k), again using a sample of size N-1, but this time we use the point (observation) i=2 for testing purposes. We will repeat this procedure N times. It is necessary to ensure that each point of a sample has been used for testing only once. Thus, we will get N values of Di(k) (i=1, N). After averaging, we will get:

\begin{equation}
\bar{D}(k) = \frac{\sum_{i=1}^{N} D_i(k)}{N} \label{eq:6} \tag{6}
\end{equation}

One of the graphs in Figure 5 (a, b, c) depicts the curve $\bar{D}(k)$. Using this curve we can determine the value K1 (the maximum number of terms of an approximating series), as well as the value S1 (the estimate of residual variance). Further, we can calculate the value of the multiple determination coefficient for each separate observation that was not used to estimate the regression function, average it over all observations and obtain the function R(k) for different values of k.

Using its maximum, we can determine K1, which in turn enables us to calculate that number of terms of an approximating series which gives the unbiased estimate of regression coefficients. Of course, the realization of this algorithm is more time consuming. Indeed, it requires N times the amount of time normally needed to perform of regression analysis. But this algorithm is a valuable tool when the sample size N is not very large. Besides, everything has its price. In this case the price is the time of calculation using a computer. But this is less expensive than creating an additional testing sample.

\subsection*{2.2 Algorithm 3 Suggested Method}
With this algorithm the computation time increases only two times in comparison to the usual method of regression analysis.

First, from the main sample we create a random sub sample of size, which is 80\% of the main sample. It means that we randomly select 80\% of all observations. Then we estimate two functions: Y1(k) using the full sample and Y2(k) using the sub-sample. After that, we can calculate the values of the following expression:

\begin{equation}
D2(k) = \frac{\sum_{i=1}^{N} (Y1_i(k) - Y2_i(k))^2}{N} \label{eq:7} \tag{7}
\end{equation}

The values of D2(k) will be similar to the values of D(k). Now we can use the form of the curve D2(k) to determine K1. Next, we calculate the value of the multiple determination coefficient R(K1) for k=K1 and compute the residual variance D(K1):

\begin{equation}
D(K1) = \frac{\sum_{i=1}^{N} (Y_i - Y1_i(K1))^2}{N} \label{eq:8} \tag{8}
\end{equation}

These values are unbiased estimates which should remain the same if we use another population or calculate forecasts.

It should be noted that the results achieved on the basis of these two algorithms were similar to the results obtained on the basis of estimation using a testing sample. This conclusion follows from the comparisons based both on real problems and on models with different levels of noise and different sample sizes.

If the curve D(k) follows form (b) or (c) from Figure 5, it means that the sample size is sufficient, and the results obtained on the basis of two algorithms will be similar to the results obtained using the standard programs SPSS and SAS. Finally, if the curve D(k) has the form (a) from figure 5, then it is better to increase the sample size. If it is not possible, then algorithm 1 or 2 must be used. If the sample size is not large, it is better to use algorithm 1.

When we know the value K1 (the number of terms on the approximation series) we can calculate the residuals (S)

\begin{equation}
S_i(K1) = Y_i - Y_i (K1) \label{eq:9} \tag{9}
\end{equation}

and then calculate the coefficient of correlation between the dependent variable Y: and the residual ($S_i$). If the value of the coefficient of correlation is different from zero, it shows that not all of the independent variables, or some functions of the values of the independent variables (squares, cubes, other functions, etc.), were included in the regression models (approximating series) because the sample size was not representative. The different value of the coefficient of correlation from zero can be found by using the hypothesis tests about the correlation coefficient (Norusis, p. B-18).

If the curve from Figure 5 takes the form of (b) or (c), and the coefficient of correlation between the dependent variable and the residual is high enough, then not all factors (independent variables) were selected for the analysis.

\section*{3. CONCLUSION}
It should be noted that to use regression analysis and standard programs, it is necessary first of all to determine the sample size in accordance with the nature of the problem. It can be done with the help of algorithm 1 or 2. If the sample size is not sufficient and if we use standard programs, then none of the coefficients of the regression equation correspond to a real model.

Moreover, the regression function can include variables which a real model would not have. On the contrary, variables which should be included in the real model will be excluded from the regression equation. In this case it is better to increase the number of observations in the sample. If it is not possible, then it is necessary to determine the number of terms of an approximating series (with the help of algorithms 1 and 2) and those variables which should be included in the regression function. Only after this procedure is completed will it be possible to begin the regression analysis using standard programs.

\section*{References}
Draper, N.R., and Smith, H. (1981), \textit{Applied Regression Analysis}, New York: Wiley. \\
Norusis, M.J. (1990), \textit{SPSS/PC and Statistics 4.0}, Chicago: SPSS, Inc. \\
Snee, R.D. (1977), \textit{Technometrics}, 415-428. 

\end{document}